\documentclass[12pt]{article}

\usepackage{amsfonts}
\usepackage{latexsym,amsmath,amssymb}

\newtheorem{thm}{Theorem}
\newtheorem{cor}[thm]{Corollary}

\newtheorem{prop}[thm]{Proposition}
\newtheorem{defn}[thm]{Definition}

\newtheorem{rem}[thm]{Remark}

\begin{document}

\title{On some aspects of the geometry of non integrable distributions and applications}


\author{Miguel-C. Mu\~noz-Lecanda\\
Departamento de Matem\'aticas-UPC\\
C. J. Girona, 3, Edif. C-3, Campus Nord-UPC\\
E-08034-Barcelona, Spain
\footnote{Email: miguel.carlos.munoz@upc.edu}}        
\date{\today} 
\maketitle
\abstract{ We consider a regular distribution $\mathcal{D}$ in a Riemannian manifold $(M,g)$. The Levi-Civita connection on $(M,g)$ together with the orthogonal projection allow to endow the space of sections of $\mathcal{D}$ with a natural covariant derivative, the intrinsic connection. Hence we have two different covariant derivatives for sections of $\mathcal{D}$, one directly with the connection in  $(M,g)$  and the other one with this intrinsic connection. Their difference is the second fundamental form of $\mathcal{D}$ and we prove it is a significant tool to characterize the involutive and the totally geodesic distributions and to give a natural formulation of the equation of motion for mechanical systems with constraints.
The two connections also give two different notions of curvature, curvature tensors and sectional curvatures, which are compared in this paper with the use of the second fundamental form. }

\bigskip
\noindent\textit{Keywords}: non-involutive distributions; Riemannian curvature; second fundamental form; totally geodesic distributions.

\bigskip
\noindent \textit{Mathematics Subject Classification (2010)}: 53C21, 58A30, 70G45, 70Q05

\section{Introduction}

Geometry of surfaces in the ordinary three dimensional space was developed by Gauss at the beginning of nineteenth century introducing the notions of first and second fundamental forms and curvature. This development was enlarged with the study of submanifolds of a Riemannian manifold. A detailed study of this subject is contained, for example, in \cite{HICKS-1967,dCA-1992,LEE-1997}. But 
there are other aspects that need to be studied as we are going to see.

In fact, any regular surface $S$ in $\mathbb{R}^3$ can be locally described in the form  $f=0$, where $f:U\rightarrow \mathbb{R}$ is a differentiable function in an open set $U\subset\mathbb{R}^3$. Taking the differential form $\alpha=\mathrm{d}\, f$, then $S\cap U$ is a solution of the  Pfaff system defined by $\alpha=0$. Let $\mathcal{D}=\mathrm{span}\{X_1,X_2\}$, with $X_i$ vector fields in $U$, be the incident distribution to $\alpha$, $\mathrm{i}(X_i)\alpha=0$, $i=1,2$.  Then $S\cap U$ is an integral manifold of the distribution $\mathcal{D}$: for every $p\in S\cap U$ we have that $ T _p S=\mathrm{span}\{X_1(p),X_2(p)\}$. Hence any surface in $\mathbb{R}^3$ can be given locally as an integral manifold of an integrable Pfaff system or its equivalent rank two distribution.

In the same way, if $N$ is a regular submanifold of $\mathbb{R}^m$, with $\mathrm{dim}\,N=n$, then $N$ is locally defined  by $F^{-1}(0)$ where  $F:U\rightarrow \mathbb{R}^{m-n}$, $F=\{f^1,\ldots ,f^{m-n}\}$, is a differentiable function. As in the previous paragraph,  $N\cap U$ is a solution of the integrable Pfaff system defined by $\Omega=\mathrm{span}\{\mathrm{d}\, f^1,\ldots,\mathrm{d}\, f^{m-n}\}$. We can also consider the distribution $\mathcal{D}$ annihilating $\Omega$; then $N\cap U$  is an integral manifold of a regular distribution $\mathcal{D}$ with rank $n$. In this case both $\Omega$ and $\mathcal{D}$ are integrable, see \cite{BULE-2005, LEE-2013, NISCH-1990} for details. In summary, any regular submanifold of $\mathbb{R}^m$ is locally an integral manifold of a Pfaff system or its incident, or dual, regular distribution.

For any submanifold $N$ of a differentiable manifold $M$ the situation is the same. So, in some sense, the geometric study of regular distributions or Pfaff systems is a natural generalization of the geometric study of submanifolds. The case of integrable Pfaff systems, or integrable distributions, corresponds to the study of foliations. See \cite{REIN-1983} for a systematic study of this subject.

The aim of this paper is to develop suitable tools to extend and study the usual notions associated to curvature of submanifolds of $\mathbb{R}^m$ to the case of regular distributions, not necessarily integrable, in a Riemannian manifold. 

Given a distribution $\mathcal{D}$ in a Riemannian manifold $(M,g)$, the Levi-Civita connection $\nabla$ on the manifold $M$ and the orthogonal projection onto the distribution are used to obtain the so called intrinsic connection on the distribution, called $\nabla^{\mathcal{D}}$. The difference between both connections when acting on sections of the distribution is a tensor field on the distribution, called the second fundamental form as in the classical study of submanifolds of a Riemannian manifold. Essentially the second fundamental form acting on two sections of the distribution is the covariant derivative with the connection $\nabla$ composed with the projection onto the  subspace orthogonal to the distribution. 

If the distribution $\mathcal{D}$ is non integrable, then this tensor field is not symmetrical and its decomposition in symmetric and skew-symmetric parts is related with different geometric properties of the distribution. This second fundamental form is our basic tool to obtain our results and we study and prove different properties of both symmetric and skew-symmetric components.

The skew-symmetric part is related to the non integrability of $\mathcal{D}$.
In fact we prove that the skew-symmetric part is zero, that is the second fundamental form is symmetric, if and only if the distribution is involutive. 

On the other hand, the symmetric part of the second fundamental form is related with the curvature of curves in the distribution and using this symmetric part we characterize in several forms the so called totally geodesic distributions. Some of these characterizations are new and for the others the proofs are considerably shortened.

The connection $\nabla$ and the intrinsic connection $\nabla^{\mathcal{D}}$ also
define  two different notions of curvature on the space of sections of the distribution, as in the case of a submanifold in a Riemannian manifold. We compare the corresponding Riemann curvature tensors and sectional curvatures obtaining several results, all of them related with the classical ones for submanifolds but with new terms coming from the non involutivity of the distribution. These new terms include the antisymmetric part of the fundamental form.

Not only is this study interesting  from the geometric viewpoint, but its applications spread to several other fields, for example to dynamical systems with controls. These controls are modeled as vector fields on the phase space of the system and one needs to know as much as possible about the specific  properties of the distribution these control vector fields span, especially if one is interested, for example, in the controllability or reachability for those systems, where the property of being totally geodesic for the control distribution is fundamental. The same situation arises for mechanical systems with controls. A detailed study  of this subject can be found in \cite{BULE-2005} and \cite{NISCH-1990} where some geometric tools are developed and applied.

A further example is the study of mechanical systems with constraints. In several cases these constraints are modeled by distributions in a Riemannian manifold, whose geometry it is necessary to know. There are two different approaches to these kinds of  systems, the nonholonomic and the vakonomic systems, the bibliography on that is extensive. Some recent references are \cite{ARKONE-1993, KUPOL-2001,BULE-2005, deLE-2012, BLMMM2012}.

Subriemannian geometry is another topic where this approach may contribute to clarifying the geometric background. We have not studied the application to this topic and it is an open question for the future.

There are some previous approaches to the study of distributions with this approach, for example \cite{GRAY-1967, REIN-1977, REIN-1983, GIL-1983}, but in all of them the definition of the second fundamental form contains only the symmetric part. Hence, from the very beginning the non integrability of the distribution is not taken into account as part of the geometrical properties of the associated tensor fields. When appropriate, we will comment on the relationships between their results, and others, with ours.

The paper is organized as follows: The second section introduces the elements to study, notation and terminology related to the problem. The third is devoted to the construction of a natural connection associated with a distribution $\mathcal{D}$ in a Riemannian manifold.

In the fourth section, we introduce the second fundamental form for a distribution and study its properties, in particular its decomposition into symmetric and skew-symmetric parts, both on the distribution $\mathcal{D}$ and on the Pfaff system associated to $\mathcal{D}$. We prove that a distribution is integrable if and only if the skew-symmetric component of the second fundamental form is identically zero.

As an application we describe the equation of motion for a mechanical system in a Riemannian manifold, both without constraints and constrained by holonomic and nonholonomic constraints, according to ideas coming at least in time from \cite{KUPOL-2001}.

The fifth section is dedicated to the study of the curvature of curves in the distribution and the notion of totally geodesic distribution, this property being characterized in different forms. In particular we prove that one distribution is totally geodesic if and only if the symmetric part of its second fundamental form is identically zero. 

As a consequence we obtain a new and shorter proof for the relation of being totally geodesic with the so called symmetric product,  \cite{BULE-2005,LEW-1997, BALE-2012}.

In the sixth section the two different notions of Riemann curvature tensor are introduced, the first one with the Levi-Civita connection $\nabla$ of $(M,g)$ and the second one with the intrinsic connection $\nabla^{\mathcal{D}}$ of the distribution. The relation of both curvature tensors is obtained by clarifying the role of the non involutivity of the distribution and the importance of the second fundamental form. We conclude this section with the study of the sectional curvature of a plane section of the distribution both with the distribution intrinsic connection and with the Levi-Civita connection of the manifold.

Finally, the seventh section is devoted to comments, open problems and future perspectives.

\section{Regular distributions in a Riemannian manifold and associated elements}

Let $(M,g)$ be a smooth Riemannian manifold, $\mathrm{dim}\,M=m$, and $\nabla$  the Levi-Civita connection associated to the Riemannian metric $g$. We denote by $\mathfrak{X}(M)$ the $\mathcal{C}^{\infty}(M)$-module of vector fields on $M$ and by $\nabla_XY$ the covariant derivative of $Y$ with respect to $X$ if  $X,Y\in\mathfrak{X}(M)$.

 If $\gamma:I\subset\mathbb{R}\to M$ is  a parametrized smooth curve, then $\mathfrak{X}(M,\gamma)$ is the $\mathcal{C}^{\infty}(I)$-module of vector fields in $M$ along $\gamma$. If $V\in\mathfrak{X}(M,\gamma)$, then $\nabla_{\dot{\gamma}}V\in\mathfrak{X}(M,\gamma)$ denotes the covariant derivative of $V$ along $\gamma$. In particular we put $\nabla_{\dot{\gamma}}\dot{\gamma}\in\mathfrak{X}(M,\gamma)$ for what is usually called the acceleration of $\gamma$.

  Let $\mathcal{D}\subseteq TM$ be a regular distribution on $M$; that is, a subbundle of the tangent bundle $TM$, with constant rank $n$. If it is necessary for the local study, we will assume that $\mathcal{D}=\mathrm{span}\{X_1,\ldots, X_n\}$, with $X_i\in\mathfrak{X}(M)$ linearly independent. In some cases we assume that these local generators are orthogonal or orthonormal. If $p\in M$, we call $\mathcal{D}_p$ the fibre of $\mathcal{D}$ at the point $p\in M$; that is the elements in $T_p M$ belonging to $\mathcal{D}$.

This distribution $\mathcal{D}$ is the object under study throughout the paper and we seek to describe its geometry using the Riemannian notions from the manifold $M$, which is the metric, as well as the Levi-Civita connection and the curvature, and the similar elements induced on $\mathcal{D}$.

Henceforth we fix the manifold $(M, g,\nabla)$ and the distribution $\mathcal{D}$ with the above conditions. All the manifolds and mappings will be regarded as being of $\mathcal{C}^{\infty}$ class.

We use these other associated elements:
\begin{enumerate}
\item $\Omega\subseteq T^*M$ is the \textbf{incident distribution} to  $\mathcal{D}$, or the annihilator  of $\mathcal{D}$. If $\mathcal{D}=\mathrm{span}\{X_1,\ldots, X_n\}$, then we put $\Omega=\mathrm{span}\{\alpha^1,\ldots,\alpha^{m-n} \}$ satisfying $\mathrm{i}(X_i)\alpha^j=0$ for every $i\in\{1,\dots,n\},j\in\{1,\ldots,m-n\}$. The distribution $\Omega$ is a subbundle of $ T ^*M$, the cotangent bundle to $M$, with constant rank $m-n$.

The sections of $\Omega$ form a Pfaff system equivalent to the distribution $\mathcal{D}$.

\item $\mathcal{D}^{\bot}\subseteq TM$ is the \textbf{orthogonal  distribution} to $\mathcal{D}$. It is a regular distribution with rank $m-n$ and it is given locally by $\mathcal{D}^{\bot}=\mathrm{span}\{Z_1,\ldots,Z_{m-n}\}$ with $g(X_i,Z_j)=0$, for every $i,j$.

\item $\Omega^{\bot}\subseteq T ^*M$ is the \textbf{incident distribution} to $\mathcal{D}^{\bot}$. When necessary, we locally write $\Omega^{\bot}=\mathrm{span}\{\beta^1,\ldots,\beta^n\}$, with $\mathrm{i}(Z_j)\beta^i=0$, for every $i,j$.

\item We say that a vector field $X\in\mathfrak{X}(M)$ is a \textbf{section} of $\mathcal{D}$ if $X(p)\in\mathcal{D}_p$ for every $p\in M$. We denote by $\Gamma(\mathcal{D})$ the $\mathcal{C}^{\infty}(M)$-module of sections of $\mathcal{D}$.

The same notation will be used for the $\mathcal{C}^{\infty}(M)$-modules of sections of the other subbundles we have defined; that is, we put  $\Gamma(\mathcal{D}^{\bot})$, $\Gamma(\Omega)$ and $\Gamma(\Omega^{\bot})$ for the sections of those bundles.

 \item
 The orthogonal decomposition $TM=\mathcal{D}\oplus\mathcal{D}^{\bot}$ gives us \textbf{natural projections}
 $$\pi^{\mathcal{D}}:TM\to\mathcal{D},\qquad\pi^{\mathcal{D}^{\bot}}:TM\to\mathcal{D}^{\bot},$$
  directly extended to the module of sections; we will use the same notation for the projections in both cases.

\item
The \textbf{musical isomorphisms}, \textsl{flat} and \textsl{sharp}, associated to $g$; that is, $g^{\flat}: T \,M\to T ^*\, M$ and its inverse $g^{\sharp}: T ^*\,M\to T \, M$, act on the above subbundles in the following way:
$$
g^{\flat}:\mathcal{D}\to\Omega^{\bot},\,\,\quad g^{\sharp}:\Omega\to\mathcal{D}^{\bot}.
$$
and they give rise to an identification between $\Omega^{\bot}$ and $\mathcal{D}^*$, the dual of $\mathcal{D}$, and between $\Omega$ and $(\mathcal{D}^{\bot})^*$.

These isomorphisms can be extended in the natural way to the sets of the corresponding sections; we will use the same notation for these extensions.
\end{enumerate}

\section{Connection induced on $\mathcal{D}$ by $(M,g,\nabla)$ }

There are two different ways to construct a covariant derivative on the set of sections of $\mathcal{D}$. One is derived from the Levi-Civita connection on the manifold $(M,g)$ and the orthogonal projection on the distribution and the other considering $\mathcal{D}$ as a Riemannian algebroid. At the end both connection are the same. We quickly review the construction of both and their specific interesting properties. 

\subsection{The intrinsic connection induced by $(M,g,\nabla)$}

The distribution $\mathcal{D}$, as a subbundle of the tangent bundle $TM$ of the manifold $M$, inherits a metric tensor field from the original $g$ in $M$. This metric will be denoted by $g^{\mathcal{D}}$ and is Riemannian in the sense that it is symmetric, non degenerate and positive definite.

Using the Levi-Civita connection $\nabla$ in $(M,g)$ and the projection $\pi^{\mathcal{D}}:TM\to\mathcal{D}$, one can define a covariant derivative between sections of  $\mathcal{D}$ in the following way:
\begin{eqnarray*}
\nabla^{\mathcal{D}}:&\Gamma(\mathcal{D})\times\Gamma(\mathcal{D})&\longrightarrow\,\,\,\,\,\Gamma(\mathcal{D})\\
 &\,\,(X,Y)&\longmapsto\pi^{\mathcal{D}}(\nabla_XY).
\end{eqnarray*}

This covariant derivative is usually called the \textbf{intrinsic connection} of the distribution $\mathcal{D}$ in the Riemannian manifold $(M,g)$.  Apart from being $\mathbb{R}$-linear in both entries, it satisfies the usual properties of a connection in a vector bundle: For every $X,Y\in\Gamma(\mathcal{D})$, $f\in\mathcal{C}^{\infty}(M)$ we have

\bigskip
\hspace{16mm}a)  $\nabla^{\mathcal{D}}_{fX}Y=f\nabla^{\mathcal{D}}_XY$,\hspace{15mm}
b)  $\nabla^{\mathcal{D}}_{X}(fY)= (L_X\,f)Y+f\nabla^{\mathcal{D}}_XY$.

\bigskip
\noindent where $L_X$ is the usual Lie derivative with respect to the vector field $X$.

Observe that the connection $\nabla^{\mathcal{D}}$ is not torsionfree with respect to the Lie bracket in the manifold $M$. In fact, for $X,Y\in\Gamma(\mathcal{D})$, the torsion of $\nabla^{\mathcal{D}}$ is given by:
\begin{equation}
 \label{eqn:tor}
T(X,Y)=\nabla^{\mathcal{D}}_X Y- \nabla^{\mathcal{D}}_Y X-[X,Y]=-\pi^{\mathcal{D}^{\bot}}([X,Y])
\end{equation}
which is zero if and only if the distribution $\mathcal{D}$ is integrable. See \cite{LEW-1998-2} for a detailed study. 

The connection $\nabla^{\mathcal{D}}$ extends to sections of the dual bundle of  $\mathcal{D}$ and to all the tensor fields on the distribution. With this in mind, it is easy to prove that $\nabla^{\mathcal{D}}$  is Riemannian with respect to $g^{\mathcal{D}}$, that is $\nabla^{\mathcal{D}}_Xg^{\mathcal{D}}=0$ for every $X\in\Gamma(\mathcal{D})$.  But, as we have seen, it is not torsionless with respect to the ordinary Lie bracket unless the distribution $\mathcal{D}$ is involutive.

The initial idea of this connection goes back at least to \cite{SYN-1926, SYN-1928} where the problem of the geodesics in a non integrable distribution was studied for the first time together with the description in geometric terms of the nonholonomic mechanical systems. Another approach can be seen in \cite{VRANCE-1926}. A modern view with some other applications can be found in \cite{LEW-1998, LEW-1998-2, BULE-2005}.

\subsection{The Levi-Civita connection}

As a further approach consider now the following operation on $\Gamma(\mathcal{D})$: For $X,Y\in\Gamma(\mathcal{D})$, define
$$
[X ,Y ]^{\mathcal{D}}=\pi^{\mathcal{D}}([X ,Y ]).
$$
This bracket $[\, ,\, ]^{\mathcal{D}}$ has the same properties as the ordinary Lie bracket except for the Jacobi identity which is satisfied if and only if the distribution $\mathcal{D}$ is involutive.

Now $(\mathcal{D}, g^{\mathcal{D}},[\, ,\,]^{\mathcal{D}})$ is a so-called skew-symmetric Riemannian algebroid and there exists the so called \textbf{Levi-Civita connection} denoted by $\bar{\nabla}^{\mathcal{D}}$ and defined by the classical Koszul formula: For every $X,Y,Z\in\Gamma(\mathcal{D})$  put
 \begin{eqnarray*}
2g^{\mathcal{D}}(\bar{\nabla}^{\mathcal{D}}_XY,Z)&=&L_X(g^{\mathcal{D}}(Y,Z))+L_Y(g^{\mathcal{D}}(X,Z))-L_Z(g^{\mathcal{D}}(X,Y))+\\
&+&g^{\mathcal{D}}(X, [Z ,Y ]^{\mathcal{D}}) +g^{\mathcal{D}}(Y,[Z,X ]^{\mathcal{D}})-g^{\mathcal{D}}(Z,[Y,X ]^{\mathcal{D}}),
\end{eqnarray*}
and $\bar{\nabla}^{\mathcal{D}}_XY$ is well defined because $g^{\mathcal{D}}$ is non degenerate.

This connection on $(\mathcal{D},g^{\mathcal{D}} )$ is the only one which is Riemannian  and torsionless with respect to the bracket  $[\, ,\, ]^{\mathcal{D}}$; that is,  the connection $\bar{\nabla}^{\mathcal{D}}$ is determined by the following properties:
\begin{center}
\begin{enumerate}
\item $[X,Y ]^{\mathcal{D}}=\bar{\nabla}^{\mathcal{D}}_XY-\bar{\nabla}^{\mathcal{D}}_YX$
\item $L_X(g^{\mathcal{D}}(Y,Z))=g^{\mathcal{D}}(\bar{\nabla}^{\mathcal{D}}_XY,Z)+g^{\mathcal{D}}(Y,\bar{\nabla}^{\mathcal{D}}_XZ)$.
\end{enumerate}
\end{center}

See \cite{BLMMM2012} for more details on this approach.

\vspace{1mm}
\begin{rem}{Comparison between both connections: $\nabla^{\mathcal{D}}=\bar{\nabla}^{\mathcal{D}}$}
\end{rem}
As the Levi-Civita connection is unique with these last two properties, to prove that both connections, $\nabla^{\mathcal{D}}$ and  $\bar{\nabla}^{\mathcal{D}}$, are the same it is enough to show that the first  one fulfills the two properties characterizing the second one  which is a straighforward exercise considering the expression (\ref{eqn:tor}). Hence $\nabla^{\mathcal{D}}=\bar{\nabla}^{\mathcal{D}}$.

Henceforth, we will denote it by $\nabla^{\mathcal{D}}$.

\section{Second fundamental form. Shape operator}

Given $X,Y\in\Gamma(\mathcal{D})\subset\mathfrak{X}(M)$, we can derivate with two different connections: if we consider them as elements in $\mathfrak{X}(M)$ with $\nabla$, but as elements of $\Gamma(\mathcal{D})$ with $\nabla^{\mathcal{D}}$; the difference of both derivatives gives us some properties of $\mathcal{D}$.

\begin{defn}

 The \textbf{second fundamental form} of $\mathcal{D}$  is the map
\begin{eqnarray*}
B:&\Gamma(\mathcal{D})\times\Gamma(\mathcal{D})&\longrightarrow\,\,\,\,\,\,\,\,\mathfrak{X}(M)\\
 &\,\,(X,Y)&\longmapsto\nabla_XY-\nabla^{\mathcal{D}}_XY=\pi^{\mathcal{D}^{\bot}}(\nabla_XY).
\end{eqnarray*}
\end{defn}
\vspace{-1mm}This map satisfies the following properties:
\begin{enumerate}
\item
$B$ is $\mathcal{C}^{\infty}(M)$-bilinear and takes values in $\Gamma(\mathcal{D}^{\bot})\subset\mathfrak{X}(M)$.
\item By definition, we have that
\vspace{-3mm}$$
\nabla_XY=\nabla^{\mathcal{D}}_XY+B(X,Y),
$$
which is known as the \textbf{Gauss formula}.

\item
Let $Z_1,\ldots,Z_{m-n}$ be a local orthonormal basis of $\Gamma(\mathcal{D}^{\bot})$, then:
\vspace{-3mm}
$$
B(X,Y)=\sum_{j=1}^{m-n}g(B(X,Y),Z_j)Z_j
$$

\vspace{-4mm}
or including it in the Gauss formula:

\vspace{-4mm}
$$
\nabla_XY=\nabla^{\mathcal{D}}_XY+\sum_{i=1}^{m-n}g(B(X,Y),Z_j)Z_j.
$$
\end{enumerate}
This is  the generalization of the known Gauss expression for surfaces in $\mathbb{R}^3$, or in general for a submanifold of a Riemannian manifold; see \cite{LEE-1997}. The proof is direct taking into account that $Z_1,\ldots,Z_{m-n}$ is an orthonormal basis of $\Gamma(\mathcal{D}^{\bot})$.

\bigskip
\begin{rem}
\end{rem}
\begin{enumerate}

\item The second fundamental form was introduced by Gauss for a surface in the ordinary three dimensional space and is used to describe, in a more general context, the geometry of a submanifold of a Riemannian or a pseudo-Riemannian manifold; see for example \cite{HICKS-1967, CHEEB-1975, SPI-1979, ONEILL-1983, LEE-1997}. However, to the best of our knowledge the oldest definition for a general distribution is given in 
\cite{REIN-1977} where, as we pointed out in the introduction, only the symmetric part is considered. Other references for its application to the study of the geometry of distributions  in Riemannian geometry have been given in the Introduction. We consider that our definition takes into account all the geometry of the distribution and is the classical one in the particular case of integrable distributions or foliations, \cite{HICKS-1967}. 

In \cite{KUPOL-2001}, the tensor field $B$ is used as an auxiliary element to write in a geometric form the dynamical equations for a mechanical system with nonholonomic constraints, both in the nonholonomic and in the vakonomic formulations. No analysis of the properties of $B$ is included, and our definition is not exactly the same because they use $B:\mathfrak{X}(M)\times\Gamma(\mathcal{D})\longrightarrow\mathfrak{X}(M)$, that is they derivate with all the vector fields of the manifold $M$ as in a classical connection in a vector bundle. In the sequel, we will comment on some of the applications studied in \cite{KUPOL-2001}.

\item
In the recent paper \cite{PRIN-2016} there is an interpretation of the antisymmetric part of the second fundamental form as torsion of a connection  but in a more general context.

\item
In reference \cite{BEFAR-2006}  the aim is to use ideas coming from nonholonomic geometry; that is, non-integrable distributions, to study different structures in foliations; that is, integrable distributions.  A description of the second fundamental form in the case of nonholonomic distribution is included in the first chapter in a similar way to the above.

\end{enumerate}

\subsection{Shape operator}
\bigskip
The above last expression of $B(X,Y)$ leads us to the following definition:
\begin{defn}
Given $X,Y\in\Gamma(\mathcal{D})$, $Z\in\Gamma(\mathcal{D}^{\bot})$, the  \textbf{shape operator} of the distribution $\mathcal{D}$ is the map
$$
S(X,Y,Z)= g(B(X,Y),Z).
$$
\end{defn}

\medskip
Other useful expressions of $S(X,Y,Z)$ are the following:
\begin{eqnarray}\label{B_Z}
S(X,Y,Z) &=&g(B(X,Y),Z)= g(\nabla_XY-\nabla^{\mathcal{D}}_XY,Z)\nonumber\\
&=&g(\nabla_XY,Z)\label{B_Z-1}\\
&=&L_X(g(Y,Z))-g(Y,\nabla_XZ)\nonumber\\
&=&-g(\nabla_XZ,Y)\label{B_Z-2}\\
&=&-g(\nabla_XZ,Y)+g(\nabla_ZX,Y)-g(\nabla_ZX,Y)\nonumber\\
&=&g(\nabla_ZX-\nabla_XZ,Y)-g(\nabla_ZX,Y)=g([Z,X],Y)-g(\nabla_ZX,Y)\nonumber\\
&=&g(L_ZX-\nabla_ZX,Y)=g((L_Z-\nabla_Z)X,Y)\label{B_Z-3}.
\end{eqnarray}

\bigskip
The following theorem collects the necessary properties of the map $S$.
\medskip
\begin{thm}

The mapping $S:\Gamma(\mathcal{D})\times\Gamma(\mathcal{D})\times\Gamma(\mathcal{D}^{\bot})\longrightarrow\mathcal{C}^{\infty}(M)$,
defined by:
$$
S(X,Y,Z)=g(B(X,Y),Z)=-g(\nabla_XZ,Y)=g(\nabla_XY,Z),
$$
satisfies the following properties:
\begin{enumerate}
  \item $S$ is $\mathcal{C}^{\infty}(M)$-linear in its three entries.
  \item $S(X,Y,Z)+S(Y,X,Z)=-(L_Z\,g)(X,Y)$.
  \item $S(X,Y,Z)-S(Y,X,Z)=-(\mathrm{d}(\mathrm{i}_Z\,g))(X,Y)$.
\end{enumerate}
\end{thm}
\textit{Proof}:
\begin{enumerate}
  \item It is direct from the definition and the different expressions  (\ref{B_Z-1}), (\ref{B_Z-2}) and (\ref{B_Z-3}) above.

\item From the above expressions of $S(X,Y,Z)$ and the properties of the Levi-Civita connection $\nabla$ on $M$, we have:
\begin{eqnarray*}
S(X,Y,Z)+S(Y,X,Z)
&=&g(L_Z X-\nabla_Z X,Y)+g(L_Z Y-\nabla_Z Y,X)\\
&=&g(L_Z X,Y)+g(L_Z Y,X)-(g(\nabla_Z X,Y)+g(\nabla_Z Y,X))\\
&=&g(L_Z X,Y)+g(L_Z Y,X)-L_Z(g(X,Y))\\
&=&g(L_Z X,Y)+g(L_Z Y,X)-((L_Zg)(X,Y)\\
&\quad&\qquad+\,\, g(L_Z X,Y)+g(L_Z Y,X))\\
&=&-(L_Z\, g)(X,Y).
\end{eqnarray*}

\item  In the same way as in the previous item:
 \begin{eqnarray*}
  S(X,Y,Z)-S(Y,X,Z)
  &=&g(\nabla_XY,Z)-g(\nabla_YX,Z)\\
  &=&g([X,Y],Z)\\
  &=&-(\mathrm{d}(\mathrm{i}_Z\,g))(X,Y)
\end{eqnarray*}
 because $(\mathrm{i}_Z\,g)(X)=g(Z,X)=0$, $(\mathrm{i}_Z\,g)(Y)=g(Z,Y)=0$.
\end{enumerate}
\vspace{-5mm}\hspace{140mm}$\blacksquare$

\bigskip
\begin{rem}:
Due to the  $\mathcal{C}^{\infty}(M)$-linear dependency of $S$ on its three entries, that is $S$ is a tensor field, if $p\in M$, $u,v\in\mathcal{D}_p$, $w\in\mathcal{D}_p^{\bot}$, then we can define $S(u,v,w)$ by using extensions of the tangent vectors $u,v,w$ as sections of $\mathcal{D}$ and $\mathcal{D}_p^{\bot}$ respectively,  and the value of  $S(u,v,w)$ is independent of the extensions we use.
\end{rem}

\begin{defn}
For a fixed $Z\in\Gamma(\mathcal{D}^{\bot})$, we denote by $B_Z$ the tensor field defined by:

\begin{eqnarray*}
B_Z:&\Gamma(\mathcal{D})\times\Gamma(\mathcal{D})&\longrightarrow\,\,\,\,\,\,\mathcal{C}^{\infty}(M)\\
  &\,\,(X,Y)&\longmapsto S(X,Y,Z).
 \end{eqnarray*}
\end{defn}

\subsection{Symmetries of $B$ and $B_Z$}

From the last two items in the above theorem, we can study the decomposition of $B$ and $S$ into their symmetric and antisymmetric components in the following way:

\begin{defn}
The \textbf{symmetric and antisymmetric components} of $B$, denoted by $B^s$ and $B^a$, are given by
$$
B^s(X,Y)= \frac{1}{2}(B(X,Y)+B(Y,X)),\,\,\,\,\,
B^a(X,Y)= \frac{1}{2}(B(X,Y)-B(Y,X)),
$$
and we have that $B=B^s+B^a$.

\end{defn}
\begin{defn}
For every $Z\in\Gamma(\mathcal{D}^{\bot})$, the \textbf{symmetric and antisymmetric components }of $B_Z$, denoted by $B_Z^s$ and $B_Z^a$, are given by
$$
B_Z^s(X,Y)= \frac{1}{2}(B_Z(X,Y)+B_Z(Y,X)),\,\,\,\,\,
B_Z^a(X,Y)= \frac{1}{2}(B_Z(X,Y)-B_Z(Y,X)),
$$
and we have that $B_Z=B_Z^s+B_Z^a$.
\end{defn}

From the definitions it is easy to verify that:

\begin{cor}\mbox

\begin{enumerate}
  \item $B=0$ if and only if $B_Z=0$ for every $Z\in\Gamma(\mathcal{D}^{\bot})$.
  \item $B^s=0$ if and only if $B_Z^s=0$ for every $Z\in\Gamma(\mathcal{D}^{\bot})$.
  \item $B^a=0$ if and only if $B_Z^a=0$ for every $Z\in\Gamma(\mathcal{D}^{\bot})$.
\end{enumerate}
\end{cor}

As we will see in the sequel, these components are necessary to characterize some of the properties of the distribution $\mathcal{D}$. First we obtain nice expressions for $B_Z^s$ and $B_Z^a$.

\begin{prop}\label{SYMSKEW}

For every $Z\in\Gamma(\mathcal{D}^{\bot})$, the components of $B_Z$, are given by:
   $$
B_Z^s=-\frac{1}{2}(L_Z\, g), \,\,\,\,\,\,
B_Z^a=-\frac{1}{2}(\mathrm{d}\,(\mathrm{i}_Z\,g)).
$$
 Hence $B_Z=-\frac{1}{2}\left(L_Z\, g+\mathrm{d}\,(\mathrm{i}_Z\,g)\right)$.

\end{prop}

\textit{Proof}: They are simple consequences of the properties of $\mathbf{S}$ we have proven in the above Theorem.

\hspace{140mm}$\blacksquare$

\bigskip
Now we have the following result

\begin{thm}

The distribution $\mathcal{D}$ is involutive if and only if for every $Z\in\Gamma(\mathcal{D}^{\bot})$ the tensor field $B_Z$ is symmetric; that is, $\mathrm{d}(\mathrm{i}_Z\,g)$ is null on the sections of $\mathcal{D}$.
\end{thm}

\textit{Proof}: For every $Z\in\Gamma(\mathcal{D}^{\bot})$, if $X,Y\in\Gamma(\mathcal{D})$, we have:
    $$
  (\mathrm{d}(\mathrm{i}_Z\,g))(X,Y)=L_X(\mathrm{i}_Z\,g(Y))-L_Y(\mathrm{i}_Z\,g(X))-\mathrm{i}_Z\,g([X,Y])=
  $$
  $$
  =L_X(g(Z,Y))-L_Y(g(Z,X))-g(Z,[X,Y])=-g(Z,[X,Y]).
  $$
  then $\mathcal{D}$ is involutive if and only if  $(\mathrm{d}(\mathrm{i}_Z\,g))(X,Y)=0$, because $Z\in\Gamma(\mathcal{D}^{\bot})$, hence if and only if  $B_Z^a=0$.

\hspace{140mm}$\blacksquare$

\begin{cor} The distribution $\mathcal{D}$ is involutive if and only if its second fundamental form  $B$ is symmetric.
\end{cor}

\medskip
If we consider the Pfaff system $\Omega$ incident to $\mathcal{D}$, then the involutivity of the distribution $\mathcal{D}$ is equivalent to saying that $\mathrm{d}\alpha(X,Y)=0$ for every $\alpha\in\Omega$ and $X,Y\in\mathcal{D}$, see \cite{CON-2001}.  In this way the last statement is no more than the Riemannian interpretation of the Frobenius theorem because  $\Omega=\{\mathrm{i}_Z\,g; Z\in\Gamma(\mathcal{D}^{\bot})\}$.

\begin{rem}\mbox\\
\end{rem}
\begin{enumerate}
\item
Weingarten map was clasically defined as: For every $Z\in\Gamma(\mathcal{D}^{\bot})$, consider the $\mathcal{C}^{\infty}(M)$-linear map:
$$
W_Z:\Gamma(\mathcal{D})\longrightarrow\mathfrak{X}(M),\quad W_Z(X)=-\nabla_X Z
$$
which satisfies $W_Z=L_Z-\nabla_Z$. Then $S(X,Y,Z)= g(Y,W_Z(X))$.
\item
If we take normalized sections $Z\in\Gamma(\mathcal{D}^{\bot})$, then $W_Z$ takes values in $\Gamma(\mathcal{D})$.
\item
For oriented surfaces $S\subset\mathbb{R}^3$, the usual way is to take the unitary normal vector field along the surface; thus, $W_Z$ takes values in $\mathfrak{X}(S)$.
\end{enumerate}

\subsection{Dual viewpoint}

The sharp isomorphism $g^{\sharp}$ gives an identification between $\mathcal{D}^{\bot}$ and $\Omega$, and this identification goes to the modules of sections of these bundles, hence we can dualize the shape operator $S$ in the previous paragraph changing its action on $\mathcal{D}^{\bot}$ by the corresponding action on $\Omega$. 

In this way we can define a tensor field 
$S^*:\Gamma(\mathcal{D})\times\Gamma(\mathcal{D})\times\Gamma(\Omega)\longrightarrow\mathcal{C}^{\infty}(M)$, by putting
$$
S^*(X,Y,\alpha)=-((\nabla_X\alpha)(Y)=\alpha(\nabla_XY)=S(X,Y,Z),
$$
where $Z=g^{\sharp}(\alpha)=\mathrm{i}_{\alpha}g^{-1}
$.  Observe that the  the third expression comes from
\vspace{-3mm}$$0=\nabla_X(\alpha(Y))=(\nabla_X\alpha)(Y)+\alpha(\nabla_XY),$$
and  the $\mathcal{C}^{\infty}(M)$-linearity of $S^*$ is clear from the last expression of $S^*$.
 
The decomposition of $S^*$ in symmetric and skew-symmetric components is given by:
\begin{eqnarray*}
S^*(X,Y,\alpha)+S^*(X,Y,\alpha) & =  S(X,Y,Z)+S(Y,X,Z) & =  -(L_Z\,g)(X,Y) \\
S^*(X,Y,\alpha)-S^*(X,Y,\alpha) & =  S(X,Y,Z)-S(Y,X,Z) & =  -(\mathrm{d}(\mathrm{i}_Z\,g))(X,Y),
\end{eqnarray*}

\noindent where, once again, the relation between $Z$ and $\alpha$ is given by $\alpha= \mathrm{i}_Zg$.

Hence for a given $\alpha\in\Gamma(\Omega)$, denoting by $S^*_{\alpha}=S^*(.,.,\alpha)$, we have that:
$$
S^*_{\alpha}=-\frac{1}{2}\left(L_{\mathrm{i}_{\alpha}g^{-1}}g+\mathrm{d}\,\alpha\right)=
-\frac{1}{2}\left(L_{g^{\sharp}(\alpha)}g+\mathrm{d}\,\alpha\right)
$$
as the decomposition of $S^*_{\alpha}$ into symmetric and antisymmetric components. This is equivalent to the expression previously obtained for $B_Z$, with $Z=\mathrm{i}_{\alpha}g^{-1}\in\Gamma(\mathcal{D}^{\bot})$, in Proposition \ref{SYMSKEW}.

Following this last comment, the tensor field $
S^*_{\alpha}$ is symmetric for every $\alpha$, if and only if the distribution $\mathcal{D}$ is integrable.

\subsection{One simple application to classical mechanics}

Let $(M,G)$ be a Riemannian manifold. A mechanical system on $M$ is given by a vector field  $F\in\mathfrak{X}(M)$, the field of force. A curve $\gamma:I\subset\mathbb{R}\to M$ is a trajectory of the system defined by $(M,g,F)$ if it is a solution to the so-called dynamical equation or Newton  equation
$$
\nabla_{\dot{\gamma}}\dot{\gamma}=F\circ\gamma.
$$
In the case that the force $F$ depends on the velocities, then it is a vector field  on $M$ along the natural projection of the tangent bundle $\tau_M:TM\to M$,  $F\in\mathfrak{X}(M,\tau_M)$, then in the dynamical equation the term $F\circ\gamma$ must be changed to $F\circ\dot{\gamma}$.

If the system has holonomic constraints, it is bounded to move on a submanifold $N\subset M$, so we must assume that there exists a new force, the constraint or reaction force $R$, and therefore make use of the \textbf{ d'Alembert principle}: the constraint force is orthogonal to the constraint submanifold $N$. Hence the dynamical equation for the trajectories of the system  $\gamma:I\subset\mathbb{R}\to N\subset M$ is:
$$
\nabla_{\dot{\gamma}}\dot{\gamma}=F\circ\gamma+R
$$
where both the curve $\gamma$ and the reaction force $R$ are unknowns. Observe that in this case $\gamma$ is a curve in the submanifold $N$.

The decomposition of this equation into tangent  and normal components to $N$ gives us two equations:
$$
\pi^N(\nabla_{\dot{\gamma}}\dot{\gamma})=\pi^N(F\circ\gamma),\qquad
\pi^{N^{\bot}}(\nabla_{\dot{\gamma}}\dot{\gamma})=\pi^{N^{\bot}}(F\circ\gamma)+R
$$
where  $\pi^N$ and $\pi^{N^{\bot}}$ are the orthogonal projections associated to the decomposition $T_pM=T_pN+(T_pN)^{\bot}$, for every $p\in N$.

Note that the second equation is no more than
$$
B({\dot{\gamma}},\dot{\gamma})=\nabla_{\dot{\gamma}}\dot{\gamma}-
\pi^N(\nabla_{\dot{\gamma}}\dot{\gamma})=\pi^{N^{\bot}}(F\circ\gamma)+R\circ\dot{\gamma}
$$
where $B$ is the second fundamental form of the submanifold $N\subset M$; compare with \cite{LEE-1997}. In this case $B$ is symmetric because $N$ is a submanifold. It is one leave of the foliation defined by an integrable distribution.

The first equation allows us to obtain the trajectory $\gamma$ and the second one to obtain the reaction force but only \textbf{along} $\dot{\gamma}$. In fact we have:
$$
R(t)=
\pi^{N^{\bot}}(\nabla_{\dot{\gamma}}\dot{\gamma})(t)-\pi^{N^{\bot}}(F\circ\gamma)(t)=B({\dot{\gamma}},\dot{\gamma})(t)-\pi^{N^{\bot}}(F\circ\gamma)(t).
$$

If the mechanical system $(M,g,F)$  is constrained in a nonholonomic way, then  there exists a distribution $\mathcal{D}\subset TM$ such that the trajectories of the system $\gamma$ must satisfy the condition $\dot{\gamma}(t)\in\mathcal{D}_{\gamma(t)}$ for every value of the parameter. Once again we must admit the \textbf{d'Alembert principle}, and thus we assume that there exists a constraint force $R$, forcing the system to satisfy the constraints and this force is orthogonal to $\mathcal{D}$. Thus, the dynamical equation consists in looking for curves  $\gamma:I\subset\mathbb{R}\to M$ such that $\dot{\gamma}(t)\in\mathcal{D}_{\gamma(t)}$, for every $t\in I$, and satisfying the equation:
$$
\nabla_{\dot{\gamma}}\dot{\gamma}=F\circ\gamma+R.
$$
Once again by orthogonal decomposition we have two different equations:
$$
\pi^{\mathcal{D}}(\nabla_{\dot{\gamma}}\dot{\gamma})=\pi^{\mathcal{D}}(F\circ\gamma),\qquad
\pi^{{\mathcal{D}}^{\bot}}(\nabla_{\dot{\gamma}}\dot{\gamma})=B({\dot{\gamma}},\dot{\gamma})
=\pi^{{\mathcal{D}}^{\bot}}(F\circ\gamma)+R
$$
where $\pi^{\mathcal{D}}$ and $\pi^{\mathcal{D}^{\bot}}$ are the projections from the tangent bundle to $\mathcal{D}$ and ${\mathcal{D}}^{\bot}$, respectively.

The first equation allows us to obtain the curve $\gamma$ and the second one to obtain the constraint force \textbf{along} $\dot{\gamma}$. Compare with \cite{KUPOL-2001, GLAKO-2004A, GLAKO-2004B,KOBOL-2004} where different, but similar, formulations are used and the second fundamental form of the non-integrable distribution  $\mathcal{D}$ is an ingredient of the dynamical equation\cite{KUPOL-2001}.

Note that in both cases, the holonomic and the nonholonomic, the dynamical equation depends only on the symmetric part of the second fundamental form. Hence they have the same expression in both cases, the integrable and the non-integrable one.

\section{Curvature of curves in the distribution}

Let $\gamma:I\subseteq\mathbb{R}\mapsto M$ be a smooth curve parametrized by the arc. The \textbf{curvature} of  $\gamma$ in $M$ is defined as $\mathbf{k}(\gamma)=||\nabla_{\dot{\gamma}}\dot{\gamma}||$.

\bigskip
In the case that $\dot{\gamma}(t)\in\mathcal{D}_{\gamma(t)}$ for all $t\in I$, we say that $\gamma$ is a \textbf{curve of the distribution} $\mathcal{D}$, and we define:
\begin{enumerate}
\item
 The \textbf{geodesic curvature} of $\gamma$ as $\mathbf{k}^{\mathcal{D}}(\gamma)=||\nabla^{\mathcal{D}}_{\dot{\gamma}}\dot{\gamma}||$.
\item
 The \textbf{normal curvature} of $\gamma$ as $\mathbf{k}^{\mathcal{D}^{\bot}}(\gamma)=||\pi^{\mathcal{D}^{\bot}}(\nabla_{\dot{\gamma}}\dot{\gamma})||=||B(\dot{\gamma},\dot{\gamma})||$.
 \end{enumerate}
 Both are functions of the parameter of the curve.

 \bigskip
 Following the definition of the second fundamental form and the Gauss formula, we have that if $\dot{\gamma}(t)\in\mathcal{D}_{\gamma(t)}$ for all $t\in I$,
 $$
 (\mathbf{k}(\gamma))^2=(\mathbf{k}^{\mathcal{D}}(\gamma))^2+(\mathbf{k}^{\mathcal{D}^{\bot}}(\gamma))^2=(\mathbf{k}^{\mathcal{D}}(\gamma))^2+||B(\dot{\gamma},\dot{\gamma})||^2.
 $$

 \begin{rem}
The normal curvature only depends on the symmetric part of the second fundamental form $B$, because we need to calculate only $B(\dot{\gamma},\dot{\gamma})$.
 \end{rem}

\subsection{Geodesics in $\mathcal{D}$}

Let $\gamma:I\subseteq\mathbb{R}\mapsto M$ be a smooth curve.

\begin{defn}

\begin{enumerate}
\item The curve $\gamma$ is \textbf{$\nabla$-geodesic} if $\nabla_{\dot{\gamma}}\dot{\gamma}=0$
\item The curve $\gamma$ is \textbf{$\nabla^{\mathcal{D}}$-geodesic} if $\dot{\gamma}(t)\in\mathcal{D}_{\gamma(t)}$, for all $t\in I$, and $\nabla^{\mathcal{D}}_{\dot{\gamma}}\dot{\gamma}=0$.
\end{enumerate}
\end{defn}

In this paragraph we are interested in the comparison between the $\nabla^{\mathcal{D}}$-geodesics and the $\nabla$-geodesics when they have initial condition in $\mathcal{D}$.

As usual, the geodesic curves are solutions to a second order ordinary differential equation whose solutions with initial condition points in  $TM$ for the $\nabla$-geodesics or points in $\mathcal{D}\subset TM$ for the $\nabla^{\mathcal{D}}$-geodesics. The existence of solutions for such equations is a consequence of their regularity and the appropriate theorem for ordinary differential equations.

The first relation between both geodesics is the following result:

\begin{prop}

Let $\gamma$ be a smooth curve in the distribution, that is $\dot{\gamma}(t)\in\mathcal{D}_{\gamma(t)}$, for all $t\in I$. Then it is  a $\nabla$-geodesic if and only if it is a $\nabla^{\mathcal{D}}$-geodesic and $B(\dot{\gamma},\dot{\gamma})=0$.
\end{prop}

\textsl{Proof}: By the definition of the second fundamental form, we have:
$$
\nabla_{\dot{\gamma}}\dot{\gamma}=\nabla^{\mathcal{D}}_{\dot{\gamma}}\dot{\gamma}+B(\dot{\gamma},\dot{\gamma})
$$
because $\dot{\gamma}(t)\in\mathcal{D}_{\gamma(t)}$. Since the last two summands are orthogonal, the conclusion in immediate.

\hspace{140mm}$\blacksquare$

\subsection{Totally geodesic distributions}

For a submanifold of a Riemannian manifold, an interesting property is for it to be totally geodesic; that is every geodesic with initial condition in the submanifold, lies locally in the submanifold, see  \cite{LEE-1997,CHEEB-1975,SPI-1979, ONEILL-1983}. For distributions we can state the same problem. In the integrable case, it corresponds to the study of the leaves of the foliation defined by the distribution. In the non-integrable situation, apart from the geometric problems, the property is specially interesting in several other fields; for example in the study of controllability of dynamical systems; see for instance \cite{BULE-2005} where they use the name geodesically invariant instead of totally geodesic. This is also part of the wide field of subriemannian geometry.

\begin{defn}

The distribution $\mathcal{D}$ is  \textbf{totally geodesic} if every $\nabla$-geodesic with initial condition in $\mathcal{D}$ is contained in $\mathcal{D}$.
\end{defn}

\begin{thm}\label{totgeod}

For a distribution $\mathcal{D}$, the following conditions are equivalent:
\begin{enumerate}
\item  $\mathcal{D}$ is totally geodesic.
\item Every $\nabla^{\mathcal{D}}$-geodesic in $\mathcal{D}$  is $\nabla$-geodesic.
\item The symmetric part of the second fundamental form is identically zero.
\item If $X,Y\in\Gamma(\mathcal{D})$ then $\nabla_XY+\nabla_YX\in\Gamma(\mathcal{D})$.
\item If $X\in\Gamma(\mathcal{D})$ then $\nabla_XX\in\Gamma(\mathcal{D})$.
\end{enumerate}
\end{thm}

\textsl{Proof}:

$\mathbf{1\Longleftrightarrow 3}$

Given $u_p\in\mathcal{D}_p$, let $\gamma$ be the $\nabla$-geodesic with initial condition $u_p\in\mathcal{D}_p$. If $\mathcal{D}$ is totally geodesic, then $\dot{\gamma}(t)\in\mathcal{D}_{\gamma(t)}$. Hence we can use the Gauss formula,
$$
0=\nabla_{\dot{\gamma}}\dot{\gamma}=\nabla^{\mathcal{D}}_{\dot{\gamma}}\dot{\gamma}+B(\dot{\gamma},\dot{\gamma})
$$
and both summands being orthogonal, we have
$$
\nabla^{\mathcal{D}}_{\dot{\gamma}}\dot{\gamma}=0,\qquad B(\dot{\gamma},\dot{\gamma})=0.
$$
In particular $B(u_p,u_p)=0$, but $u_p$ is every point in $\mathcal{D}$, hence $B^s=0$.

Conversely, assume that $B^s=0$. Given $u_p\in\mathcal{D}$, consider the following curves: $\gamma$ the $\nabla$-geodesic and $\sigma$ the $\nabla^{\mathcal{D}}$-geodesic, both with initial condition $u_p\in\mathcal{D}_p$. We have that $\dot{\sigma}(t)\in\mathcal{D}_{\sigma(t)}$, hence by the Gauss formula, being $B^s=0$, we have
$$
\nabla_{\dot{\sigma}}\dot{\sigma}=\nabla^{\mathcal{D}}_{\dot{\sigma}}\dot{\sigma}+B(\dot{\sigma},\dot{\sigma})=\nabla^{\mathcal{D}}_{\dot{\sigma}}\dot{\sigma}=0.
$$
Then $\gamma=\sigma$ by unicity of solutions. From this we have that every $\nabla$-geodesic beginning in $\mathcal{D}$ is contained in $\mathcal{D}$ as we sought to prove.

$\mathbf{2\Longleftrightarrow 3}$

If $B^s=0$, for every curve $\gamma$ in ${\mathcal{D}}$ we have $\nabla_{\dot{\gamma}}\dot{\gamma}=\nabla^{\mathcal{D}}_{\dot{\gamma}}\dot{\gamma}$, hence every $\nabla^{\mathcal{D}}$-geodesic in ${\mathcal{D}}$ is $\nabla$-geodesic.

Conversely, if every curve $\gamma$ being $\nabla^{\mathcal{D}}$-geodesic in $\mathcal{D}$ is $\nabla$-geodesic, then  $B(\dot{\gamma},\dot{\gamma})=0$ for all of them. However, every point in ${\mathcal{D}}$ can be taken as initial condition, so $B^s=0$.

$\mathbf{3\Longleftrightarrow 4}$

Given $Z\in\Gamma(\mathcal{D}^{\bot})$, take $X,Y\in\Gamma(\mathcal{D})$, then the symmetric component of $B_Z$ satisfies:
$$
B_Z^s(X,Y)=B_Z(X,Y)+B_Z(Y,X)=g(\nabla_XY+\nabla_YX,Z).
$$
Hence $B^s_Z=0$ if and only if  $\nabla_XY+\nabla_YX$ is orthogonal to $Z$ for every $X,Y\in\Gamma(\mathcal{D})$. As this is true for every $Z\in\Gamma(\mathcal{D}^{\bot})$, we have the result.

$\mathbf{4\Longleftrightarrow 5}$

If $\nabla_XY+\nabla_YX\in\Gamma(\mathcal{D})$ for every $X,Y\in\Gamma(\mathcal{D})$, then trivially $\nabla_XX\in\Gamma(\mathcal{D})$ taking $Y=X$.

Conversely, assume that $\nabla_XX\in\Gamma(\mathcal{D})$ for every $X\in\Gamma(\mathcal{D})$. Then given $X,Y\in\Gamma(\mathcal{D})$, we have that $\nabla_{(X+Y)}(X+Y)\in\Gamma(\mathcal{D})$, but
$$
\nabla_{(X+Y)}(X+Y)=\nabla_XX+\nabla_XY+\nabla_YX+\nabla_YY
$$
then $\nabla_XY+\nabla_YX\in\Gamma(\mathcal{D})$ as we wanted.

And this finishes the proof.

\hspace{140mm}$\blacksquare$

\bigskip
The expression $\nabla_XY+\nabla_YX$ for sections in the distribution $\mathcal{D}$, or in manifold $M$ with connection, is known as the \textbf{symmetric product}. It was introduced in \cite{CRO-1981} and has been deeply studied in  \cite{LEW-1998-2} and \cite{BULE-2005} including some of its applications. According to this nomenclature, the fourth condition  of the above theorem can be stated as in the next corollary:

\begin{cor}\cite{LEW-1997, LEW-1998}

The distribution $\mathcal{D}$ is totally geodesic if and only if $\Gamma(\mathcal{D})$ is closed under the symmetric product.
\end{cor}

\begin{rem}
\end{rem}
\begin{enumerate}
\item
The proof given above is considerably shorter than the previous proof in the references
 \cite{BALE-2012, BULE-2005} for distributions, where as we stated above they refer to it as geodesically invariant instead of our totally geodesic,  and in references \cite{dCA-1992,{CHEEB-1975}} for submanifolds.
\item Observe the importance of the second fundamental form for submanifolds and its substitution by the symmetric part in the case of non-involutive distributions.
The reason for this substitution is that every property related to the curvature of curves in the distribution $\mathcal{D}$ depends only on the symmetric part of the second fundamental form.
But when considering the sectional curvatures, as in the next section, we will show they depend on the full second fundamental form including its skew-symmetric component.
\item In \cite{GRAY-1967, GIL-1983}, the equivalence   
$\mathbf{1\Longleftrightarrow 3}$ is proved in the cases they study. In these papers the second fundamental form is defined as the symmetric part of our definition.
Given the distribution $\mathcal{D}$ in a Riemannian manifold, $M$, they consider the tensor field $P:TM\to TM$ defined by $P(u)=u$ if $u\in\mathcal{D}$, and $P(v)=-v$ if $v\in\mathcal{D}^{\bot}$. Then they prove that:
\begin{enumerate}
\item The distribution $\mathcal{D}$ is integrable if and only if $(\nabla_XP)Y=(\nabla_YP)X$ for all $X,Y\in \Gamma(\mathcal{D})$.
\item The distribution $\mathcal{D}$ is totally geodesic if and only if  $(\nabla_XP)X=0$ for every $X\in\Gamma(\mathcal{D})$.
\end{enumerate}
and this second part is equivalent to our result $\mathbf{1\Longleftrightarrow 3}$. No relation with the symmetric product is given.

\item
Condition 3 in the above theorem is not equivalent to stating that the vector fields $Z\in\Gamma(\mathcal{D}^{\bot})$ are Killing vector fields for the Riemannian metric $g$, because condition $L_Zg=0$ is true only when acting on sections of ${\mathcal{D}}$, and not for every vector field in the manifold $M$.

However, if one considers an integrable distribution such that its orthogonal distribution is made by Killing vector fields, then $B^S_Z=-(1/2)L_Zg=0$. Hence the leaves are totally geodesic.

The same is true for a hypersurface in a Riemannian manifold: if there is a Killing vector field that is normal to the submanifold, then the surface is totally geodesic. See \cite{SPI-1979} for some comments on this situation.
\end{enumerate}

\subsection{The case $n=m-1$: curvature of curves}

If the rank of  $\mathcal{D}$  is $n=m-1$, then we can take a local unitary normal section $N\in\Gamma(\mathcal{D}^{\bot})$, $g(N,N)=1$, and obtain a scalar fundamental form, $\mathbf{b}:\Gamma(\mathcal{D})\times\Gamma(\mathcal{D})\to\mathcal{C}^{\infty}(M)$, in the following way:
$$
\mathbf{b}(X,Y)=S(X,Y,N)=g(B(X,Y),N)=g(\nabla_XY,N)=-g(Y,\nabla_XN).
$$
Observe that the unitary normal section can be chosen arbitrarily as $\pm N$, and this option changes the sign of $\mathbf{b}$ in the same way.

Relating $\mathbf{b}$ with the classical Weingarten map $W:\Gamma(\mathcal{D})\rightarrow\Gamma(\mathcal{D})$, we have that
$$
W(X)=-\nabla_XN=(L_N-\nabla_N)(X),
$$ 
because $\nabla$ is the Levi-Civita connection in $(M,g)$. Hence $W=L_N-\nabla_N$ as in the case of an integrable distribution, but in our case $W$ is not symmetric with respect to $g$. See \cite{HICKS-1967} for details in the integrable case. Thus we have:
$$
\mathbf{b}(X,Y)=-g(Y,\nabla_XN)=g(Y,W(X))=g(Y,(L_N-\nabla_N)X)
$$

Using the decomposition of the fundamental form into symmetric and skew-symmetric components, we obtain for  $\mathbf{b}$:
\begin{eqnarray*}
\mathbf{b}(X,Y)+\mathbf{b}(X,Y)&=&-(L_Ng)(X,Y)\\
\mathbf{b}(X,Y)-\mathbf{b}(X,Y)&=&-(\mathrm{d}\,\mathrm{i}_N g)(X,Y),
\end{eqnarray*}
as in the general case. That is:
$$
\mathbf{b}^s=-\frac{1}{2} L_Ng,\,\,\,\,\,\, \mathbf{b}^a=-\frac{1}{2}\mathrm{d}\,\mathrm{i}_N g.
$$
and $\mathbf{b}^s$ and $\mathbf{b}^a$ are the symmetric and skew-symmetric parts of $\mathbf{b}$. In the case of an integrable distribution, $B^a$ is zero, hence $\mathbf{b}^a$ is zero and $\mathbf{b}$ is symmetric. Recall that $\mathbf{b}$, and hence $\mathbf{b}^s$ and $\mathbf{b}^a$, are tensor fields acting only on $\Gamma(\mathcal{D})$, not on $\mathfrak{X}(M)$, the set of all the vector fields on the manifold $M$.

As we know, the curvature of curves in $\mathcal{D}$ depends only on the symmetric part of  the second fundamental form, then only on $\mathbf{b}^s$ in this case. Hence, to study the curvature of curves in the distribution we can work in the non integrable case, using only the symmetric part $\mathbf{b}^s$, as in the case of an integrable distribution, where we describe the curvature of the leaves of the associated foliation, see  \cite{REIN-1983}. Thus, we obtain the principal curvatures and principal directions of the distribution $\mathcal{D}$ using $\mathbf{b}^s$ and the corresponding associated symmetric endomorphism by contraction with $g^{-1}$. The difference with the classical theory of integrable distributions is that in our case we have no points or curves contained in a submanifold but only curves with tangent vector in the distribution.

Obviously the distribution $\mathcal{D}$ is integrable if and only if $\mathbf{b}^a=0$.

In the following section we will see that for the sectional curvatures  the skew-symmetric part is necessary.

\section{Curvature tensors.  Gauss theorem}

To derive the elements of $\Gamma(\mathcal{D})$ there are two connections acting on, $\nabla$ and  $\nabla^{\mathcal{D}}$. Then we have two different notions of curvature, endomorphism of curvature and curvature tensor corresponding to those two connections. For curvature definitions and notation we follow \cite{LEE-1997}.

Given $X_1,X_2,X_3\in\Gamma(\mathcal{D})$, the Riemann \textbf{curvature endomorphism}, $\mathbf{R}$, is the tensor field defined by
$$
\mathbf{R}(X_1,X_2)X_3=\nabla_{X_{1}}\nabla_{X_{2}}X_{3}-  \nabla_{X_{2}}\nabla_{X_{1}}X_{3}- \nabla_{[X_{1},X_{2}]}X_{3}
$$
$$
\mathbf{R}^{\mathcal{D}}(X_1,X_2)X_3=\nabla_{X_{1}}^{\mathcal{D}}\nabla_{X_{2}}^{\mathcal{D}}X_{3}-  \nabla_{X_{2}}^{\mathcal{D}}\nabla_{X_{1}}^{\mathcal{D}}X_{3}- \nabla_{[X_{1},X_{2}]^{\mathcal{D}}}^{\mathcal{D}}X_{3},
$$
with the connection $\nabla$ of $M$, acting on the sections of  ${\mathcal{D}}$, and with the intrinsic connection $\nabla^{\mathcal{D}}$ of the distribution $\mathcal{D}$ respectively. The first one refers to the curvature endomorphism as vector fields in the manifold $M$ and the second one as sections of the distribution $\mathcal{D}$.

On the other hand, given $X_1,X_2,X_3,X_4\in\Gamma(\mathcal{D})$, the Riemann \textbf{curvature tensor}, $\mathrm{K}$, is defined by
\begin{eqnarray*}
\mathrm{K}(X_1,X_2,X_3,X_4)&=&g(\mathbf{R}(X_1,X_2)X_3,X_4)\\
&=&g(\nabla_{X_{1}}\nabla_{X_{2}}X_{3}-  \nabla_{X_{2}}\nabla_{X_{1}}X_{3}- \nabla_{[X_{1},X_{2}]}X_{3},X_{4})
\end{eqnarray*}
\begin{eqnarray*}
\mathrm{K}^{\mathcal{D}}(X_1,X_2,X_3,X_4)&=&g(\mathbf{R}^{\mathcal{D}}(X_1,X_2)X_3,X_4)\\
&=&g(\nabla_{X_{1}}^{\mathcal{D}}\nabla^{\mathcal{D}}_{X_{2}}X_{3}-  \nabla^{\mathcal{D}}_{X_{2}}\nabla^{\mathcal{D}}_{X_{1}}X_{3}- \nabla^{\mathcal{D}}_{[X_{1},X_{2}]^{\mathcal{D}}}X_{3},X_{4})
\end{eqnarray*}
with respect to both connections.

In this paragraph we study the relations between both tensor fields in the two different connections.

For the curvature endomorphism we have the following result:

\begin{thm}\label{CURVENDOM}

\begin{enumerate}
\item
Given $X_1,X_2,X_3\in\Gamma(\mathcal{D})$, we have
\begin{eqnarray*}
\mathbf{R}^{\mathcal{D}}(X_1,X_2)X_3&=&\pi^{\mathcal{D}}(\mathbf{R}(X_1,X_2)X_3)\\
&\quad&\quad -(\pi^{\mathcal{D}}\circ\nabla_{X_{1}})(B(X_{2},X_{3}))
+(\pi^{\mathcal{D}}\circ\nabla_{X_{2}})(B(X_{1},X_{3}))\\
&\quad&\quad +\pi^{\mathcal{D}}\left(\nabla_{[X_{1},X_{2}]^{\mathcal{D}^{\bot}}}X_{3}\right).
\end{eqnarray*}
\item
If $Z_{j}$, $j=1,\ldots,m-n$, is a local orthonormal basis of $\Gamma(\mathcal{D}^{\bot})$, then:
\begin{eqnarray*}
\mathbf{R}^{\mathcal{D}}(X_1,X_2)X_3&=&\pi^{\mathcal{D}}(\mathbf{R}(X_1,X_2)X_3)\\
&\quad&\quad -\sum_{j}B_{Z_{j}}(X_{1},X_{3})\pi^{\mathcal{D}}(\nabla_{X_{2}}Z_{j})\\
&\quad&\quad +\sum_{j}B_{Z_{j}}(X_{2},X_{3})\pi^{\mathcal{D}}(\nabla_{X_{1}}Z_{j})
\\
&\quad&\quad +\pi^{\mathcal{D}}\left(\nabla_{[X_{1},X_{2}]^{\mathcal{D}^{\bot}}}X_{3}\right)
\end{eqnarray*}
where $[X_{1},X_{2}]^{\mathcal{D}^{\bot}}=\pi^{\mathcal{D}^{\bot}}([X_{1},X_{2}])$ is the natural bracket in the distribution $\mathcal{D}^{\bot}$.
\end{enumerate}
\end{thm}

\textsl{Proof}: By direct calculation, we have
\begin{eqnarray*}
\nabla_{X_{1}}^{\mathcal{D}}\nabla^{\mathcal{D}}_{X_{2}}X_{3}-  \nabla^{\mathcal{D}}_{X_{2}}\nabla^{\mathcal{D}}_{X_{1}}X_{3}&-& \nabla^{\mathcal{D}}_{[X_{1},X_{2}]^{\mathcal{D}}}X_{3}\\[1mm]
& =&\pi^{\mathcal{D}}\left(\nabla_{X_{1}}\nabla^{\mathcal{D}}_{X_{2}}X_{3}-  \nabla_{X_{2}}\nabla^{\mathcal{D}}_{X_{1}}X_{3}- \nabla_{[X_{1},X_{2}]^{\mathcal{D}}}X_{3}\right)\\
&=&(\pi^{\mathcal{D}}\circ\nabla_{X_{1}})(\nabla^{\mathcal{D}}_{X_{2}}X_{3})-
(\pi^{\mathcal{D}}\circ\nabla_{X_{2}})(\nabla^{\mathcal{D}}_{X_{1}}X_{3})\\[1mm]
& &-\pi^{\mathcal{D}}\left(\nabla_{[X_{1},X_{2}]^{\mathcal{D}}}X_{3}\right)\\[1mm]
&=&(\pi^{\mathcal{D}}\circ\nabla_{X_{1}})(\nabla^{\mathcal{D}}_{X_{2}}X_{3})-
(\pi^{\mathcal{D}}\circ\nabla_{X_{2}})(\nabla^{\mathcal{D}}_{X_{1}}X_{3})\\[1mm]
& &-\pi^{\mathcal{D}}\left(\nabla_{[X_{1},X_{2}]}X_{3}\right)+
\pi^{\mathcal{D}}\left(\nabla_{[X_{1},X_{2}]^{\mathcal{D}^{\bot}}}X_{3}\right).
\end{eqnarray*}

Using the second fundamental form in the first two terms:
\begin{eqnarray*}
(\pi^{\mathcal{D}}\circ\nabla_{X_{1}})(\nabla^{\mathcal{D}}_{X_{2}}X_{3})=
(\pi^{\mathcal{D}}\circ\nabla_{X_{1}})\left(\nabla_{X_{2}}X_{3}-B(X_{2},X_{3})\right)\\
(\pi^{\mathcal{D}}\circ\nabla_{X_{2}})(\nabla^{\mathcal{D}}_{X_{1}}X_{3})=
(\pi^{\mathcal{D}}\circ\nabla_{X_{2}})\left(\nabla_{X_{1}}X_{3}-B(X_{1},X_{3})\right)
\end{eqnarray*}
and substituting above:
\begin{eqnarray*}
\nabla_{X_{1}}^{\mathcal{D}}\nabla^{\mathcal{D}}_{X_{2}}&\!\!\!\!\!\!X_{3}&\!\!\!\!\!\!-  \nabla^{\mathcal{D}}_{X_{2}}\nabla^{\mathcal{D}}_{X_{1}}X_{3}- \nabla^{\mathcal{D}}_{[X_{1},X_{2}]^{\mathcal{D}}}X_{3}
=\\
[2mm]& &\pi^{\mathcal{D}}\left(\nabla_{X_{1}}\nabla_{X_{2}}X_{3}-  \nabla_{X_{2}}\nabla_{X_{1}}X_{3}- \nabla_{[X_{1},X_{2}]}X_{3}\right)\\[2mm]
&\quad&\quad -(\pi^{\mathcal{D}}\circ\nabla_{X_{1}})(B(X_{2},X_{3}))
+(\pi^{\mathcal{D}}\circ\nabla_{X_{2}})(B(X_{1},X_{3}))\\[2mm]
&\quad&\quad +\pi^{\mathcal{D}}\left(\nabla_{[X_{1},X_{2}]^{\mathcal{D}^{\bot}}}X_{3}\right).
\end{eqnarray*}

And this is the first part of the theorem.

For the second item, let $Z_{j}$, $j=1,\ldots,m-n$, be a local orthonormal basis of $\Gamma(\mathcal{D}^{\bot})$. Then:
\begin{eqnarray*}
\nabla_{X_{2}}(B(X_{1},X_{3}))&=&\nabla_{X_{2}}\left(\sum_{j}B_{Z_{j}}(X_{1},X_{3})Z_{j}\right)\\
&=&
\sum_{j}\nabla_{X_{2}}\left(B_{Z_{j}}(X_{1},X_{3})\right)Z_{j}+
\sum_{j}B_{Z_{j}}(X_{1},X_{3})\nabla_{X_{2}}Z_{j}.
\end{eqnarray*}
Hence:
$$
(\pi^{\mathcal{D}}\circ\nabla_{X_{2}})(B(X_{1},X_{3}))=
\sum_{j}B_{Z_{j}}(X_{1},X_{3})\pi^{\mathcal{D}}(\nabla_{X_{2}}Z_{j}),
$$
because $\pi^{\mathcal{D}}(Z_j)=0$ for every $j$. In a similar way:
$$
(\pi^{\mathcal{D}}\circ\nabla_{X_{1}})(B(X_{2},X_{3}))=
\sum_{j}B_{Z_{j}}(X_{2},X_{3})\pi^{\mathcal{D}}(\nabla_{X_{1}}Z_{j}).
$$

By substitution in the result of the first item:
\begin{eqnarray*}
\nabla_{X_{2}}^{\mathcal{D}}\nabla^{\mathcal{D}}_{X_{1}}X_{3}&-&  \nabla^{\mathcal{D}}_{X_{1}}\nabla^{\mathcal{D}}_{X_{2}}X_{3}-\nabla^{\mathcal{D}}_{[X_{1},X_{2}]^{\mathcal{D}}}X_{3}=\\[2mm]
& &\pi^{\mathcal{D}}\left(\nabla_{X_{2}}\nabla_{X_{1}}X_{3}-  \nabla_{X_{1}}\nabla_{X_{2}}X_{3}- \nabla_{[X_{1},X_{2}]}X_{3}\right)\\[2mm]
& &-\sum_{j}B_{Z_{j}}(X_{1},X_{3})\pi^{\mathcal{D}}(\nabla_{X_{2}}Z_{j})
 +\sum_{j}B_{Z_{j}}(X_{2},X_{3})\pi^{\mathcal{D}}(\nabla_{X_{1}}Z_{j})
\\
& & +\pi^{\mathcal{D}}\left(\nabla_{[X_{1},X_{2}]^{\mathcal{D}^{\bot}}}X_{3}\right)
\end{eqnarray*}
which is the expression we were looking for.

\hspace{140mm}$\blacksquare$

\begin{cor}

If $\mathbf{R}=0$, that is the ambient manifold $(M,g,\nabla)$ has zero curvature endomorphism, then:
\begin{eqnarray*}
\mathbf{R}^{\mathcal{D}}(X_1,X_2)X_3&=&
 -(\pi^{\mathcal{D}}\circ\nabla_{X_{1}})(B(X_{2},X_{3}))
+(\pi^{\mathcal{D}}\circ\nabla_{X_{2}})(B(X_{1},X_{3}))\\
&\quad&\quad +\pi^{\mathcal{D}}\left(\nabla_{[X_{1},X_{2}]^{\mathcal{D}^{\bot}}}X_{3}\right).
\end{eqnarray*}
\end{cor}

\hspace{140mm}$\blacksquare$

For the curvature tensor, we have the relation given by the following result:

\begin{thm}(\textbf{Gauss theorem})

Let $X_1,X_2,X_3,X_4\in\Gamma(\mathcal{D})$, we have that:
\begin{eqnarray}\label{GATH}
\mathrm{K}^{\mathcal{D}}(X_1,X_2,X_3,X_4)&=&\mathrm{K}(X_1,X_2,X_3,X_4)-\\
& &\; -g(B(X_{1},X_{3}),B(X_{2},X_{4}))+g(B(X_{2},X_{3}),B(X_{1},X_{4}))\nonumber\\
& &\; +g(\pi^{\mathcal{D}}(\nabla_{[X_{1},X_{2}]^{\mathcal{D}^{\bot}}}X_{3}),X_{4})\nonumber.
\end{eqnarray}
\end{thm}

\textsl{Proof}: By application of the second item of \ref{CURVENDOM}, the expression for $\mathrm{K}^{\mathcal{D}}(X_1,X_2,X_3,X_4)$ changes in the following way:
\begin{eqnarray*}
\mathrm{K}^{\mathcal{D}}(X_1,X_2,X_3,X_4)
&=&g(\nabla_{X_{2}}^{\mathcal{D}}\nabla^{\mathcal{D}}_{X_{1}}X_{3}-  \nabla^{\mathcal{D}}_{X_{1}}\nabla^{\mathcal{D}}_{X_{2}}X_{3}- \nabla^{\mathcal{D}}_{[X_{1},X_{2}]^{\mathcal{D}}}X_{3},X_{4})\\
&=&g( \pi^{\mathcal{D}}\left(\nabla_{X_{2}}\nabla_{X_{1}}X_{3}-  \nabla_{X_{1}}\nabla_{X_{2}}X_{3}- \nabla_{[X_{1},X_{2}]}X_{3}\right) ,X_{4})\\
&\quad&\quad -g(\sum_{j}B_{Z_{j}}(X_{1},X_{3})\pi^{\mathcal{D}}(\nabla_{X_{2}}Z_{j}),X_{4})\\
&\quad&\quad +g(\sum_{j}B_{Z_{j}}(X_{2},X_{3})\pi^{\mathcal{D}}(\nabla_{X_{1}}Z_{j})
 ,X_{4})\\
&\quad&\quad +g(  \pi^{\mathcal{D}}\left(\nabla_{[X_{1},X_{2}]^{\mathcal{D}^{\bot}}}X_{3}\right),X_{4})\\
&=&\mathrm{K}(X_1,X_2,X_3,X_4)-\sum_{j}B_{Z_{j}}(X_{1},X_{3})g(\pi^{\mathcal{D}}(\nabla_{X_{2}}Z_{j}),X_{4})\\
& &\; + \sum_{j}B_{Z_{j}}(X_{2},X_{3})g(\pi^{\mathcal{D}}(\nabla_{X_{1}}Z_{j}) ,X_{4})+
g(  \pi^{\mathcal{D}}\left(\nabla_{[X_{1},X_{2}]^{\mathcal{D}^{\bot}}}X_{3}\right),X_{4})\\
&=&\mathrm{K}(X_1,X_2,X_3,X_4) -\sum_{j}B_{Z_{j}}(X_{1},X_{3})B_{Z_{j}}(X_{2},X_{4})\\
& &\; + \sum_{j}B_{Z_{j}}(X_{2},X_{3})B_{Z_{j}}(X_{1},X_{4})+g(  \pi^{\mathcal{D}}\left(\nabla_{[X_{1},X_{2}]^{\mathcal{D}^{\bot}}}X_{3}\right),X_{4})\\
&=&\mathrm{K}(X_1,X_2,X_3,X_4)-g(B(X_{1},X_{3}),B(X_{2},X_{4})\\
&\quad&\quad +g(B(X_{2},X_{3}),B(X_{1},X_{4}))+g(\pi^{\mathcal{D}}\left(\nabla_{[X_{1},X_{2}]^{\mathcal{D}^{\bot}}}X_{3}\right),X_{4}).
\end{eqnarray*}
That is:
\begin{eqnarray*}
\mathrm{K}^{\mathcal{D}}(X_1,X_2,X_3,X_4)
&=&\mathrm{K}(X_1,X_2,X_3,X_4)\\
&-&g(B(X_{1},X_{3}),B(X_{2},X_{4}))+g(B(X_{2},X_{3}),B(X_{1},X_{4}))\\
&+&g(\pi^{\mathcal{D}}\left(\nabla_{[X_{1},X_{2}]^{\mathcal{D}^{\bot}}}X_{3}\right),X_{4})
\end{eqnarray*}
as we wanted.

\hspace{150mm}$\blacksquare$

\bigskip
This expression is similar to the classic one, see for example \cite{LEE-1997}, but the last term is new and is identically zero if and only if the distribution is involutive.

In \cite{GRAY-1967}, Theorem 2.4, a similar equation to \ref{GATH} is proved, but the last three terms are substituted by four terms using what is called the configuration of the distribution, a $(1,2)$ tensor field associated to the distribution,  its orthogonal complement, the natural projections and the Levi-Civita connection. 

For the particular case of distributions in $\mathbb{R}^m$, we have:

\begin{cor}

If $\mathrm{K}=0$, that is the ambient manifold $(M,g,\nabla)$ has zero curvature tensor. Then:
\begin{eqnarray*}
\mathrm{K}^{\mathcal{D}}(X_1,X_2,X_3,X_4)=
&-&g(B(X_{1},X_{3}),B(X_{2},X_{4}))+g(B(X_{2},X_{3}),B(X_{1},X_{4}))+\\
&+&g(\pi^{\mathcal{D}}\left(\nabla_{[X_{1},X_{2}]^{\mathcal{D}^{\bot}}}X_{3}\right),X_{4}).
\end{eqnarray*}
\end{cor}

\bigskip
In the case of an involutive distribution we obtain the classic results for Riemannian foliations; see \cite{REIN-1983}:

\begin{cor}\label{GAUSS-CLASICO}

Assume that the distribution $\mathcal{D}$ is involutive; then
\begin{enumerate}
\item
Given $X_1,X_2,X_3,X_4\in\Gamma(\mathcal{D})$,
\begin{eqnarray*}
\mathrm{K}^{\mathcal{D}}(X_1,X_2,X_3,X_4)
&=&\mathrm{K}(X_1,X_2,X_3,X_4)-\\
&\quad&\quad -g(B(X_{1},X_{3}),B(X_{2},X_{4}))+g(B(X_{2},X_{3}),B(X_{1},X_{4})).
\end{eqnarray*}
\item
If in addition $\mathrm{K}=0$, then
\begin{eqnarray*}
\mathrm{K}^{\mathcal{D}}(X_1,X_2,X_3,X_4)=
-g(B(X_{1},X_{3}),B(X_{2},X_{4}))+g(B(X_{2},X_{3}),B(X_{1},X_{4})).
\end{eqnarray*}
\end{enumerate}
\end{cor}
This last item is an expression of the classical \textbf{Gauss Theorema Egregium} for surfaces in $\mathbb{R}^3$.

\subsection{Sectional curvatures}

Let $X,Y\in\Gamma(\mathcal{D})$ be linearly independent at every point; the \textbf{sectional curvature} of the subspace they span is defined as
$$
\mathbf{K}(X,Y)=\frac{\mathrm{K}(X,Y,X,Y)}{g(X,X)g(Y,Y)-g(X,Y)^2}
$$
an expression that depends only on the subspace not on the specific basis we take. As in the previous cases, we can calculate it from the ambient manifold $M$; that is $\mathbf{K}$, or from the distribution $\mathcal{D}$; that is $\mathbf{K}^{\mathcal{D}}$. The expression can be simplified by taking  $X,Y$ as an orthonormal basis for the subspace.

As a consequence  of the above results, the relation between both sectional curvatures, in $M$ and in $\mathcal{D}$, are as follows:

\begin{cor}

\begin{enumerate}
\item
Given $X,Y\in\Gamma(\mathcal{D})$
\begin{eqnarray*}
\mathrm{K}^{\mathcal{D}}(X,Y,X,Y)
&=&\mathrm{K}(X,Y,X,Y)-\\
&-&g(B(X,X),B(Y,Y))+g(B(Y,X),B(X,Y))+\\
&+&g(\pi^{\mathcal{D}}\left(\nabla_{[X,Y]^{\mathcal{D}^{\bot}}}X\right),Y).
\end{eqnarray*}
\item
If in addition $\mathrm{K}=0$, we have:
\begin{eqnarray*}
\mathrm{K}^{\mathcal{D}}(X,Y,X,Y)&=&
-g(B(X,X),B(Y,Y))+g(B(Y,X),B(X,Y))+\\
&+&g(\pi^{\mathcal{D}}\left(\nabla_{[X,Y]^{c^{\bot}}}X\right),Y).
\end{eqnarray*}
\end{enumerate}
\end{cor}

Note once again that this is the classical expression,  see \cite{LEE-1997}, with the last term added. This term corresponds to the non-involutivity of $\mathcal{D}$. For the involutive situation we find the same expression as in a Riemannian foliation, that is:

\bigskip
\begin{cor}

Given $X,Y\in\Gamma(\mathcal{D})$ we have:

\begin{enumerate}
\item
If $\mathcal{D}$ is involutive, then:
\begin{eqnarray*}
\mathrm{K}^{\mathcal{D}}(X,Y,X,Y)
&=&\mathrm{K}(X,Y,X,Y)-\\
&\quad&\quad -g(B(X,X),B(Y,Y))+g(B(Y,X),B(X,Y)).
\end{eqnarray*}
\item
If in addition $\mathrm{K}=0$, then
\begin{eqnarray*}
\mathrm{K}^{\mathcal{D}}(X,Y,X,Y)=
-g(B(X,X),B(Y,Y))+g(B(Y,X),B(X,Y)).
\end{eqnarray*}
\end{enumerate}
\end{cor}
The last part of the above corollary is another expression of the classical \textbf{Gauss Theorema Egregium }: The sectional curvatures of the leaves of the foliation defined by the involutive distribution $\mathcal{D}$ do not depend on $\mathcal{D}^{\bot}$, although the second part of the last expression uses the orthogonal complement to $\mathcal{D}$.

\section{Final comments and perspectives}

Using the induced metric, the naturally induced bracket and connection, we have defined the second fundamental form for a regular distribution, whether integrable or not, in a Riemannian manifold. We have described its properties and its decomposition into symmetric and skew-symmetric components. The symmetric part has been the significant one for the study of curvature of curves in the distribution, including the characterization of the totally geodesic distributions and an adequate expression for the trajectories of a constrained mechanical system. Furthermore we have compared the curvature tensor from the point of view of the ambient manifold connection and the distribution connection. The skew-symmetric part of the second fundamental form has been important for comparing these curvature tensors and the sectional curvatures.

Now, to conclude the paper, we list some points to take into account for future work:
\begin{enumerate}
\item The vakonomic case for a mechanical system with constraints has not been included as one of the applications. The description given in \cite{KUPOL-2001} for these systems is related with specific properties of the full second fundamental form, not only with its symmetric component. In the future, we aim to describe this relation in deepth.

\item  The distribution $\mathcal{D}$ is a submanifold, a subbundle, of the tangent bundle $TM$. This bundle has a natural Riemannian metric, the Sasaki metric. Hence we can reduce our problem to the more simple problem of integrable distributions in the tangent bundle or to one of its associated leaves. At present, we are unable to state whether this approach simplifies the  description or not.

\item As stated in the introduction, subriemannian geometry is another topic where we will seek to apply the results in this paper, and we hope to obtain some results in the future.

\item The action of isometries is another problem to study. Are there any kind of invariants for classifying distributions under isometries of the Riemannian manifold $(M,g)$? The precedents are the first and second fundamental forms together with Codazzi-Mainardi conditions for surfaces in $\mathbb{R}^3$. Perhaps the approach adopted in the second item above would be more suitable for tackling this problem.

Similar problems to these last have been studied by A.  Solov$^{\prime}$ev, in (\cite{SOLOV-1982, SOLOV-1982-2}). Problems related to the classification of regular distributions on a Riemannian manifold under the action of isometries. The definitions given by Solov$^{\prime}$ev are not the same of those in this paper, not even similar to other more restricted approach as those given in (\cite{REIN-1977, REIN-1983}), and his aim is to classify the Riemannian submersions.

\end{enumerate}

\bigskip
\textbf{Acknowledgemnts}
We acknowledge the financial support of the Spanish ``Ministerio de Econom\'\i a y Competitividad"  project MTM2014-54855-P and from the Catalan Government  project 2017-SGR-932.

We thank the anonymous referee for his precise and complete report which has allowed to improve the manuscript.



\end{document}